\documentclass{amsart}

\usepackage{amsmath,amssymb}
\usepackage{graphicx}

\usepackage{array,tabularx,booktabs,enumerate}
\newcolumntype{R}{>{\raggedleft\arraybackslash}X}
\newcolumntype{C}{>{\centering\arraybackslash}X}

\usepackage{url}

\begin{document}
\allowdisplaybreaks

\title[An exact column-generation approach for the lot-type design problem]{An exact column-generation approach for the lot-type design problem -- extended abstract}

\author{Miriam Kie\ss ling \and Sascha Kurz \and J\"org Rambau}

\address{Lehrstuhl f\"ur Wirtschaftsmathematik, Universit\"at Bayreuth, 95440 Bayreuth, Germany\\
\texttt{$\{$miriam.kiessling,sascha.kurz,joerg.rambau$\}$@uni-bayreuth.de}}

\maketitle

\begin{abstract}
  We consider a fashion discounter that supplies any of its many
  branches with an integral multiple of lots whose size assortment
  structure stems from a set of many applicable lot-types.  We design
  a column generation algorithm for the optimal approximation of the
  branch and size dependent demand by a supply using a bounded number
  of lot-types.  
  
  \medskip
  
  \noindent
  \textbf{Keywords:} lot-type design, real world data, column generation
\end{abstract}

\section{Introduction}
\label{sec_introduction}

\noindent
Our business partner, a fashion discounter with more than 1500 branches,
orders all products in multiples of so-called \textit{lot-types} from
the suppliers and distributes them without any replenishing. A lot-type
specifies a number of pieces of a product for each available size, e.g., 
lot-type $(1,2,2,1)$ means two pieces of size $\text{M}$ and $\text{L}$,
one piece of size $\text{S}$ and $\text{XL}$, if the sizes are
$(\text{S},\text{M},\text{L},\text{XL})$.

We want to solve the following approximation problem: which (integral)
multiples of which (integral) lot-types should be supplied to a set of
branches in order to meet a (fractional) expected demand as closely as
possible? An ILP formulation was introduced in~\cite{p_median}, but
unfortunately for many practical instances the set of applicable
lot-types is so large that it cannot be solved directly. In this paper,
we therefore propose an exact column generation approach, which
simultaneously generates columns and cuts. For similar approaches see,
e.g., \cite{note,sim}. Unifying general
remarks can be found in \cite{0979.90092,branch_and_price_generic_scheme}. 

\section{Formal problem statement}
\label{sec_formal_problem_statement}

\textit{Data.} Let $\mathcal{B}$ be the set of branches, 
$\mathcal{S}$ be the set of sizes, and $\mathcal{M}\subset\mathbb{N}$ 
be the set of possible multiples. A \emph{lot-type} is a vector 
$(l_s)_{s\in\mathcal{S}}\in\mathbb{N}^{|\mathcal{S}|}$ satisfying 
$\min_c\le l_s\le\max_c$ for all $s\in\mathcal{S}$ and $\min_t\le
\sum_{s\in\mathcal{S}} l_s\le \max_t$.
By $\mathcal{L}$ we abbreviate the set of applicable \emph{lot-types}.
There is an upper bound $\overline{I}$ and a lower bound $\underline{I}$
given on the total supply over all branches and sizes.
Moreover, there is an upper bound $k \in \mathbb{N}$ on the number of
lot-types used. By $d_{b,s}\in\mathbb{Q}_{\ge 0}$ we denote the demand at
branch~$b$ in size~$s$.

\textit{Decisions.} Consider an assignment of a unique lot-type $l(b)
\in \mathcal{L}$ and an assignment of a unique multiplicity $m(b)\in\mathcal{M}$
to each branch~$b \in \mathcal{B}$.  These data specify that $m(b)$ lots
of lot-type $l(b)$ are to be delivered to branch~$b$.

\emph{Objective.} The goal is to find a subset $L \subseteq \mathcal{L}$ of at most $k$
lot-types and assignments $l(b) \in \mathcal{L}$ and $m(b) \in \mathcal{M}$
such that the total supply is within the bounds
$\bigl[\underline{I},\overline{I}\bigr]$, and the deviation between inventory
and demand is minimized.

\section{Modelling}
\label{sec:modelling}

\noindent
With binary assignment variables $x_{b,l,m}$ indicating whether
$l(b)=l$ and $m(b)=m$ and binary selection variables $y_l$ indicating
whether $l \in L$, we can formulate the following integer
linear program. As an abbreviation we use 
$\lvert l \rvert:=\sum\limits_{s\in\mathcal{S}}l_s$.
\begin{align}
  \label{OrderModel_Target}
  \min && \sum_{b\in\mathcal{B}}\sum_{l\in\mathcal{L}}\sum_{m\in\mathcal{M}} c_{b,l,m}\cdot x_{b,l,m}\\
  \label{OrderModel_EveryBranchOneLottype}
  s.t.  &&
  \sum_{l\in\mathcal{L}}\sum_{m\in\mathcal{M}} x_{b,l,m} &= 1 && \forall b\in\mathcal{B}\\
  \label{OrderModel_UsedLottypes}
  &&
  \sum_{l\in\mathcal{L}} y_l & \le k\\
  \label{OrderModel_Binding}
  && \sum_{m\in\mathcal{M}} x_{b,l,m} & \le y_l && \forall
  b\in\mathcal{B}, l\in\mathcal{L}
\end{align}
\begin{align}
  &&
  \label{OrderModel_Cardinality}
  \underline{I} \le
  \sum_{b\in\mathcal{B}}\sum_{l\in\mathcal{L}}\sum_{m\in\mathcal{M}}
  m \cdot \lvert l \rvert \cdot x_{b,l,m} &\le \overline{I}\\
  &&
  \label{ie_bin_x}
  x_{b,l,m} & \in\{0,1\} && \forall b\in\mathcal{B}, l\in\mathcal{L}, m\in\mathcal{M}\\
  &&
  \label{ie_bin_y}
  y_l & \in\{0,1\} && \forall l\in\mathcal{L},
\end{align}
where $c_{b,l,m}= \sum\limits_{s \in \mathcal{S}} \bigl\lvert d_{b, s} - m \cdot l_s
\bigr\rvert\ge 0$.

\section{A custom-made branch-and-price algorithm}
\label{sec_column_generation}

\noindent
The parameters $\min_c=0$, $\max_c=5$, $\min_t=12$, $\max_t=30$, and
$|\mathcal{S}|=12$ result in a set of applicable lot-types of size
1\,159\,533\,584. Thus, the number of variables and constraints of the
stated ILP formulation is, in many practical settings, very large. 
Some algorithmic observations on real-world data:
\begin{itemize}
\item The integrality gap of our ILP model is small (see \cite{p_median}).
\item Solutions generated by heuristics perform very well (see \cite{p_median}).
\item For small subsets $\bar{\mathcal{L}}\subset\mathcal{L}$ the problem can be
      solved efficiently, e.g.\ by the proposed ILP formulation.
\item A proof of optimality is wanted.
\end{itemize}

The LP relaxation of the master problem can be restricted to a manageable 
sized restricted master problem (RMP)
by the following: We consider a (small) subset
$\mathcal{L}'\subseteq\mathcal{L}$ of lot-types. For each branch
$b\in\mathcal{B}$ we consider a subset $\zeta(b)\subseteq\mathcal{L}'$
of these lot-types and for each
$l\in\zeta(b)$ we consider a subset 
$\eta(b,l)\subseteq\mathcal{M}$. To overcome the integrality gap we utilize
cover-cuts, see  Inequality~(\ref{ie_cover}), where the $\mathcal{C}_i$
are subsets of the set of lot-types.
\begin{align}
  \min && \sum_{b\in\mathcal{B}}\sum_{l\in\zeta(b)}\sum_{m\in\eta(b,l)} c_{b,l,m}\cdot x_{b,l,m}\\
  s.t.  &&
  \label{ie_one_selection_per_branch}
  \sum_{l\in\zeta(b)}\sum_{m\in\eta(b,l)} x_{b,l,m} &= 1 && \forall b\in\mathcal{B}\\
  &&
  \sum_{l\in\mathcal{L}'} -y_l & \ge -k\\
  && \!\!\!\!\!\!\!\!\!\!\!
  \overline{I}\ge\sum_{b\in\mathcal{B}}\sum_{l\in\zeta(b)}\sum_{m\in\eta(b,l)}
   m \cdot \lvert l \rvert \cdot x_{b,l,m} &\ge \underline{I}\\
  &&
  \sum_{m\in\eta(b,l)}-x_{b,l,m}+y_l & \ge 0 && \forall b\in\mathcal{B}, l\in\zeta(b)\\
  &&
  \label{ie_cover}
  \sum_{l\in\mathcal{C}_i} -y_l &\ge -(k-1) && \forall i\in\mathcal{I}\\
  &&
  x_{b,l,m} & \ge 0 && \!\!\!\!\!\!\!\!\!\!\!\!\!\!\forall b\in\mathcal{B}, l\in\zeta(b), m\in\eta(b,l) \\
  && y_l & \ge 0 && \forall l\in\mathcal{L}'.
\end{align}
This led us to the following branch-and-price algorithm:
\begin{enumerate}[(1)]
 \item Use the heuristics from \cite{p_median} to determine a 
            starting solution $(x^\star,y^\star)$.
 \item Initialize the RMP (see below),
            as follows: For each branch $b$ we compute the three
            (locally) best fitting lot-types and add them to $\zeta_b$.
            Additionally we add all lot-types used in $(x^\star,y^\star)$.
            We set $\mathcal{L}'=\cup_{b\in\mathcal{B}} \,\zeta(b)$.  For
            each branch $b\in\mathcal{B}$ and each lot-type $l\in\zeta(b)$ 
            we compute the corresponding optimal multiplicity $\hat{m}$
            and set $\eta(b,l)=\left\{\hat{m}-1,\hat{m},\hat{m}+1\right\}
            \cap\mathcal{M}$.
 \item Let $(x',y')$ be the optimal solution of RMP. If the costs are
            smaller than the costs of $(x^\star,y^\star)$,
            then we set $\bar{\mathcal{L}}=\{l\in\mathcal{L}'\mid y'_l\ge \varepsilon\}$,
            where $\varepsilon$ is a small constant, e.g., $\varepsilon=0.15$,  
            and branch on $\bar{\mathcal{L}}$, i.e.\ we perform step (5).
 \item We solve the pricing problem and possibly add lot-types from $\mathcal{L}'$
            to a $\zeta_b$, enlarge a $\eta(b,l)$, or add a new lot-type to
            $\mathcal{L}'$, i.e.\ we generate new columns and rows, go on with
            step~(3), or stop otherwise.
 \item Solve the lot-type design problem restricted to the set $\bar{\mathcal{L}}$
            of applicable lot-types and possibly update the best solution $(x^\star,y^\star)$.
            Add the cover-cut $\sum_{l\in\mathcal{C}_i} y_l \le k-1$ with
            $\mathcal{C}_i=\bar{\mathcal{L}}$ to RPM and go to step~(3).
\end{enumerate}

\section{Computational results}
\label{sec_results_column_generation}

\noindent
In this section we evaluate our proposed branch-and-price algorithm,
see Table~\ref{table_performance_CG}, using \texttt{ILOG CPLEX 12.1.0}. The key parameters
of some selected representative
problem instances, where $|\mathcal{M}|=3$, are summarized in Table~\ref{table_problem_instances}.

\begin{table}[ht]
  \begin{center}\sffamily\footnotesize
    \begin{tabularx}{\linewidth}{r*{5}{R}}
      \toprule
      Instance          &  1 &    2 &       3 &    4 &          5 \\
      \midrule
      $k$               &  3 &    5 &       5 &    4 &          5 \\
      $|\mathcal{B}|$   & 10 &   10 &    1303 & 1328 &        682 \\
      $|\mathcal{S}|$   &  4 &    4 &       4 &    7 &         12 \\
      $|\mathcal{L}|$   & 50 & 1211 &    1211 & 1290 & 1\,159\,533\,584 \\
      $[\underline{I},\overline{I}]$    & [54,66] &   [54,66] & [11\,900,12\,100] & [9702,9898] &    [15\,500,16\,200] \\
      \bottomrule
    \end{tabularx}
  \end{center}
  \caption{Key parameters for some selected problem instances.}
  \label{table_problem_instances}
\end{table}



\begin{table}[ht]
  \begin{center}\sffamily\footnotesize
    \begin{tabularx}{\linewidth}{r*{5}{R}}
      \toprule
      Instance          &  1 &   2 &       3 &    4 &       5 \\
      \midrule
      CPU~time [s]               &  1 &  1 &       2 &    4 &     937 \\
      \# variables (initial LP)  & 74 & 76 & 20\,952 & 7975 & 13\,129 \\
      \# constraints (initial LP)& 43 & 43 &    7944 & 5315 &    6642 \\
      \# variables (final LP)    & 76 & 76 & 20\,952 & 9937 & 64\,877 \\
      \# constraints (final LP)  & 44 & 43 &    7944 & 6248 & 32\,183 \\
      \# variables (complete ILP) &   1550 &  37\,541 &  3\,634\,211 & 5\,140\,650 & $2.373\cdot 10^{12}$\\ 
      \# constraints (complete ILP) &    513 &  12\,123 &  1\,212\,003 & 1\,714\,451 & 7$.908\cdot 10^{11}$ \\
      \# cover cuts              &  0 &  0 &       0 &    2 &       4 \\
      \# pricing steps           &  2 &  1 &       1 &    4 &     141 \\
      \bottomrule
    \end{tabularx}
  \end{center}
  \caption{Performance of the column generation algorithm from Section~\ref{sec_column_generation}.}
  \label{table_performance_CG}
\end{table}

The number of columns and cuts as well as the number of
branch-and-price-and-cut nodes necessary for our method to find an
optimal solution and prove its optimality is managable in all test
instances.  Instance~5 is typical for real-world data and could not be
solved statically.  Therefore, our algorithm makes optimal lot-type
design possible for industrial scale problem sizes.


\end{document}